\documentclass[12pt]{article}
\usepackage{arxiv}

\usepackage{siunitx}
\usepackage[table]{xcolor}
\definecolor{royal_blue}{rgb}{0.255, 0.412, 1}
\definecolor{medium_blue}{rgb}{0, 0, 0.804}
\usepackage{hyperref}
\hypersetup{colorlinks=true,linkcolor=medium_blue, citecolor=medium_blue}
\usepackage{graphicx}
\usepackage{bm}
\usepackage{amssymb}
\usepackage{multirow}
\usepackage{siunitx}
\usepackage{physics}
\usepackage{subcaption} 
\usepackage{float}
\usepackage{amsmath}
\usepackage{amssymb}
\usepackage{mathdots}
\usepackage{stmaryrd}
\usepackage{threeparttable}
\usepackage{makecell}
\usepackage{tabu,booktabs}
\usepackage[bottom]{footmisc}

\title{An Explicit Total Lagrangian Fragile Points Method for Finite Deformation of Hyperelastic Materials}

\author{
  Konstantinos A. ~Mountris\thanks{[mail] k.mountris@ucl.ac.uk \quad [url] https://www.mountris.org} \\
  Mechanical Engineering \\ University College London \\ London, United Kingdom \\ \texttt{k.mountris@ucl.ac.uk} \\
  \And
  Mingjing ~Li \\ School of Aeronautic Science \& Engineering \\ Beihang University \\ Beijing, China \\ \texttt{10756@buaa.edu.cn} \\
  \And
  Richard ~Schilling \\ Barts Heart Centre \\ St Bartholomew’s Hospital\\ London, United Kingdom \\
  \texttt{richard.schilling@nhs.net} \\
  \And
  Leiting ~Dong \\ School of Aeronautic Science \& Engineering \\ Beihang University \\ Beijing, China \\ \texttt{ltdong@buaa.edu.cn} \\
  \And
  Satya N. ~Atluri \\ Mechanical Engineering \\ Texas Tech University \\ Lubbock, Texas \\ \texttt{snatluri.ttu@gmail.com} \\
  \And
  Alicia ~Casals \\
  Automatic Control \& Computer Engineering \\ Universitat Politecnica de Catalunya \\ Barcelona, Spain \\ \texttt{alicia.casals@upc.edu} \\
  \And
  Helge A. ~Wurdemann \\ Mechanical Engineering \\ University College London \\ London, United Kingdom \\ \texttt{h.wurdemann@ucl.ac.uk} \\
}

\begin{document}
\maketitle

\begin{abstract}
This research explored a novel explicit total Lagrangian Fragile Points Method (FPM) for finite deformation of hyperelastic materials. In contrast to mesh-based methods, where mesh distortion may pose numerical challenges, meshless methods are more suitable for large deformation modelling since they use enriched shape functions for the approximation of displacements. However, this comes at the expense of extra computational overhead and higher-order quadrature is required to obtain accurate results. In this work, the novel meshless method FPM was used to derive an explicit total Lagrangian algorithm for finite deformation. FPM uses simple one-point integration for exact integration of the Galerkin weak form since it employs simple discontinuous polynomials as trial and test functions, leading to accurate results even with single-point quadrature. The proposed method was evaluated by comparing it with FEM in several case studies considering both the extension and compression of a hyperelastic material. It was demonstrated that FPM maintained good accuracy even for large deformations where FEM failed to converge.
\end{abstract}


\section{Introduction}
Rubber-like materials such as elastomers and soft tissue are commonly modelled under the assumption of a hyperelastic constitutive model where the stress--strain relationship is derived from a strain energy density function \cite{yeoh1997hyperelastic}. The Finite Element Method (FEM) is frequently used to simulate the elastic response of these materials when they undergo a large deformation. Although FEM is the standard numerical method for finite deformation analysis, it suffers from low convergence in conditions of large deformation and near-incompressibility \cite{hu2011meshless}. Even when non-locking elements are utilised, convergence and overall accuracy are significantly deteriorated \cite{malkus1978mixed} due to the large mesh distortion. A remedy to the mesh distortion problem is to remesh the deformed geometry during the evolution of the deformation \cite{ramalho2021novel}. However, this solution is not time efficient, especially for soft tissue simulation where results should be generated by taking into account the time restrictions of the corresponding clinical application.

Over the past several decades, meshless methods have been developed and demonstrated as being more suitable for finite deformation compared to FEM. Meshless methods either partially or fully alleviate the mesh requirement, avoiding the FEM mesh distortion problem. Examples of meshless methods include the smoothed particle hydrodynamics (SPH) \cite{monaghan2005smoothed,low2021parameter}, the meshless local Petrov-Galerkin (MLPG) \cite{atluri1998new,hu2007meshless}, and the element free Galerkin (EFG) methods \cite{belytschko1994element,wu2021interpolating}. The SPH method has been used successfully to derive a total Lagrangian formulation of SPH for large deformation simulation in 3D cardiac mechanics \cite{lluch2019breaking}. Since gradients of constant and linear functions are not correctly obtained in the standard SPH method, a corrected smoothed particle method \cite{chen1999corrective} was used to ensure the conservation of linear and angular momentum. The MLPG is another method that has been proposed as an alternative to FEM for large deformation simulations \cite{hu2011meshless,hu2007meshless}. Enriched radial basis functions with a polynomial basis function were used for constructing the MLPG trial functions. The polynomial basis function enrichment allowed the exact imposition of essential boundary conditions. Large deformation was simulated in a series of 2D benchmark problems demonstrating superior accuracy for MLPG compared to FEM.

In a similar manner to SPH and MLPG, the EFG method has been employed for solving large deformation problems. An EFG-based algorithm, the so-called Meshless Total Lagrangian Explicit Dynamics (MTLED) algorithm, has been used to a great extent for the simulation of large deformations in biomechanics \cite{joldes2017new,joldes2019suite}. The MTLED employs the EFG method to solve the total Lagrangian formulation explicitly using a central difference scheme. Initially, the method was introduced using the Moving Least Squares (MLS) for the approximation of trial functions \cite{horton2010meshless}. It is well known that MLS do not possess the Kronecker delta property, and therefore the essential boundary conditions cannot be imposed directly as in FEM. An efficient algorithm for the imposition of the exact essential boundary conditions in MTLED has been introduced in \cite{joldes2017new}. However, the application of such treatments for the imposition of essential boundary conditions can be avoided by replacing the MLS approximation functions. In \cite{mountris2020cell}, the Cell-based Maximum Entropy (CME) approximants have been proposed as an alternative to MLS. CME possesses the weak Kronecker delta property that allows applying essential boundary conditions directly and exactly as in FEM. An alternative approach is followed in \cite{bourantas2021simple}, where a regularized weight function \cite{most2005moving} is used for the construction of MLS, rendering the approximation almost interpolating. Both approaches presented accurate imposition of essential boundary conditions in several 3D benchmark problems and were cross-validated against analytical solutions and FEM \cite{mountris2020cell, bourantas2021simple}.

The aforementioned meshless methods are valuable alternatives to FEM for large deformation simulations. However, they possess drawbacks, including the requirement for complex algorithms to impose the essential boundary conditions and inexact integration due to the complexity of the approximation functions that lead to the requirement of high--order integration to improve accuracy \cite{mazzia2007comparison}. Therefore, although meshless methods can reduce dramatically the preprocessing time, their computational cost is usually significantly higher than FEM \cite{nguyen2008meshless}. Recently, a novel Fragile Points Method (FPM) that inherently avoids these limitations has been introduced by the groups of Atluri and Dong \cite{dong2019new,yang2019elementarily}. In FPM, the problem domain $\Omega$ is discretized by distributing points randomly in its interior and boundary. Following this step, $\Omega$ is partitioned in contiguous and non-overlapping subdomains that involve a single node, usually at the centroid of the subdomain. The meshless approximation is then constructed on a compact support domain that includes the points for which their subdomain shares an interface with the subdomain of the point of interest. Trial and test functions which are local, simple, polynomial, and discontinuous are employed as \cite{liszka1980finite}. Therefore, integration becomes trivial and low-order integration is sufficient to compute integrals with high accuracy in contrast to previous meshless methods. Moreover, both essential and natural boundary conditions are imposed directly as in FEM without the need for additional treatment. Nevertheless, it should be noted that deriving the standard Galerkin weak form with the FPM will lead to inconsistent and inaccurate solutions due to the discontinuity of the employed trial and test functions. To remedy this issue, interior penalty numerical flux corrections, which are commonly used in Discontinuous Galerkin methods \cite{arnold2002unified}, are employed to restore the consistency and accuracy of the FPM-derived Galerkin weak form.

Despite its recent introduction, FPM has successfully simulated different problems and provided results in agreement with FEM while reducing the requirement for a high-quality mesh. In addition, FPM was validated by solving the patch test and the Poisson equation in 1D and 2D then comparing its convergence to FEM \cite{dong2019new}. FPM demonstrated higher accuracy than FEM when mesh quality was deteriorated due to mesh distortion. Moreover, FPM has been successful in solving 2D linear elasticity \cite{yang2019elementarily} and heat conduction \cite{guan2020new} problems in 2D and 3D. Furthermore, it has been demonstrated as a promising alternative to FEM for cardiac modelling applications. Results with similar accuracy to FEM were obtained for cardiac electrophysiology simulation through the monodomain reaction-diffusion model \cite{mountris2021meshfree}. Similarly, high correlation of FPM with FEM was obtained when solved the Laplace-Dirichlet problem for the determination of cardiac muscle fiber orientation by \cite{mountris4073648meshless}. Finally, FPM possesses the ability to model crack developments in a very simple manner due to the discontinuity of the trial and test functions. It has been used successfully for modelling flexoelectric problems with crack propagation \cite{guan2021new} and damage as well as fracture of U-notched structures \cite{wang2022fragile}.

In this work, we hypothesize that FPM is a promising alternative to FEM for finite deformation simulation, especially when large mesh distortion occurs during deformation. Therefore, FPM is employed to derive the total Lagrangian formulation to simulate the large deformation of hyperelastic materials using explicit integration in time as in MTLED \cite{horton2010meshless,joldes2019suite}. Since FPM allows for direct imposition of essential and natural boundary conditions as well as accurate approximation of integrals with low-order integration rules, accurate and efficient solutions to large deformation problems are expected. The structure of the remaining paper is as follows. In section \ref{sec:theory}, the FPM total Lagrangian formulation is described, with details provided about the derivation of the explicit total Lagrangian algorithm, the formulation of the FPM trial and test functions, and the application of the interior penalty numerical flux correction. In section \ref{sec:examples}, the FPM total Lagrangian algorithm is evaluated in a series of validation case studies and the obtained solutions are compared with FEM simulations. Finally, a discussion about the findings of this work and a conclusion with future work directions are provided in sections \ref{sec:discuss} and \ref{sec:conclude}, respectively.

\section{Explicit Total Lagrangian Fragile Points Method} \label{sec:theory}
\subsection{Problem definition} \label{ss:tled}
\begin{figure}[hbt]
\captionsetup[subfigure]{justification=centering}
\begin{subfigure}[t]{0.45\textwidth}
    \includegraphics[width=\textwidth]{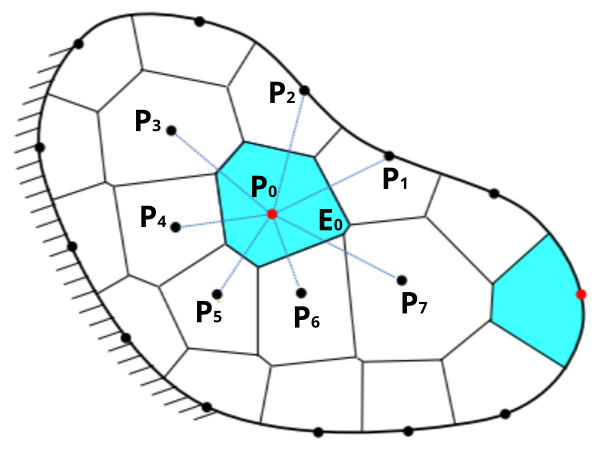}
    \caption{}
\end{subfigure}
\hfill
\begin{subfigure}[t]{0.45\textwidth}
    \includegraphics[width=\textwidth]{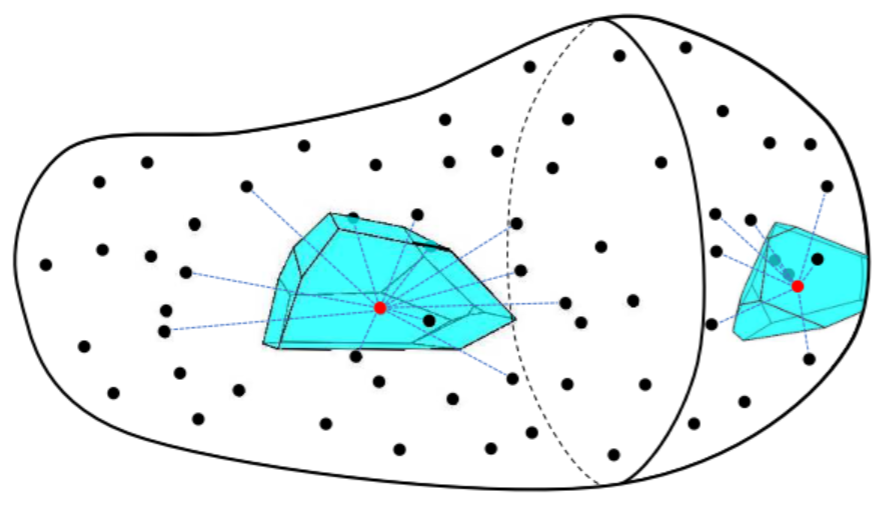}
    \caption{}
\end{subfigure}
\caption{Partitioning of domain $\Omega_0$. (a) 2D domain partitioning with points distributed inside it and on its boundary ($P \in \Omega_0 \bigcup \Gamma$). (b) 3D domain partitioning with points distributed inside it ($P \in \Omega_0$).}
\label{fig:partition}
\end{figure}

Consider the boundary value problem of an elastomer undergoing large deformation and governed by the following equations expressed in the reference configuration (total Lagrangian formulation) \cite{marsden1994mathematical}:
\begin{eqnarray} \label{eq:strong_form}
    \nabla_{0} \bm{P} + \rho_{0} \bm{b} = \rho_{0} \bm{\ddot u} &\textrm{ in } \Omega_{0} \nonumber\\
    \bm{u} = \bm{\Bar u} &\textrm{ in } \Gamma_{u} \nonumber\\
    \bm{P} \cdot \bm{n}_{0} = \bm{\Bar t} &\textrm{ in } \Gamma_{t}
\end{eqnarray}
\noindent In Equation (\ref{eq:strong_form}), $\bm{P}$ denotes the first Piola-Kirchhoff stress tensor, $\rho_0\bm{b}$ the  body forces per unit reference volume, $\ddot{\bm{u}}$ the second time derivative of the displacement, and $\bm{n}_0$ the outward unit surface normal. The subscript $0$ delineates field values expressed in respect with the reference configuration. Therefore, $\Omega_0$ is the region occupied by the body in its reference configuration while $\Gamma_u$ and $\Gamma_t$ are the parts of the boundary $\Gamma = \Gamma_u \bigcup \Gamma_t$ where boundary conditions are imposed on displacement and traction, respectively. The weak form of Equation (\ref{eq:strong_form}) is obtained by replacing $\bm{u}$ with a trial function $\bm{u}^h$ and minimizing the resulting residual through multiplication by the test function $\bm{v}$ and integration to obtain:
\begin{equation} \label{eq:weak_form}
    \int_{\Omega_0} \nabla_0 \bm{P}(\bm{u}^h) \cdot \bm{v} d\Omega + \int_{\Omega_0} \rho_0 (\bm{b} - \ddot{\bm{u}}^h) \cdot \bm{v} d\Omega = 0 \\
\end{equation}
\noindent Applying the divergence theorem and the traction boundary condition, the weak form can be rewritten as:
\begin{equation} \label{eq:weak_form_mod}
    \int_{\Omega_0} \bm{S} : \bm{E} d\Omega - \int_{\Omega_0} \rho_0 (\bm{b} - \ddot{\bm{u}}^h) \cdot \bm{v} d\Omega - \int_{\Gamma_t} \Bar{\bm{t}} \cdot \bm{v} d\Gamma = 0 \\
\end{equation}
\noindent where the first Piola-Kirchhoff stress tensor is replaced by the second Piola-Kirchhoff stress tensor ($\bm{S}$) through $\bm{P} = \bm{S} \bm{F}$ and $\bm{E}$ denotes the Green-Lagrange strain tensor obtained via $\bm{E} = \frac{1}{2} \left [ \nabla \bm{v}^T \bm{F} + \bm{F}^T \nabla \bm{v} \right]$ with $\bm{F}$ being the deformation gradient tensor. $\bm{S}$ is obtained by differentiating the material's strain energy density ($W$) by $\bm{E}$:
\begin{equation} \label{eq:strain_energy}
    \bm{S} = \frac{\partial W}{\partial \bm{E}} = \frac{\partial W}{\partial J_1} \frac{\partial J_1}{\partial \bm{E}} + \frac{\partial W}{\partial J_2} \frac{\partial J_2}{\partial \bm{E}} + \frac{\partial W}{\partial J_3} \frac{\partial J_3}{\partial \bm{E}}
\end{equation}
\noindent where $J_1$, $J_2$, $J_3$ denote the reduced invariants of the right Cauchy-Green deformation tensor. In FPM, the reference configuration domain ($\Omega_0$) is discretized by a number of randomly distributed points inside the domain and on its boundary. The domain is then further partitioned into contiguous and non-overlapping subdomains of arbitrary shapes as shown in Figure \ref{fig:partition}. Using the same type of trial function $\bm{u}^h$ and test function $\bm{v}$ in each subdomain, the discretized Galerkin weak form is obtained:
\begin{equation} \label{eq:galerkin_form}
    \sum_{E \in \Omega_0} \int_{E} \bm{S} : \bm{E} d\Omega - \sum_{E \in \Omega_0} \int_{E} \rho_0 (\bm{b} - \ddot{\bm{u}}^h) \cdot \bm{v} d\Omega - \sum_{e \in \Gamma_t} \int_{e} \Bar{\bm{t}} \cdot \bm{v} d\Gamma = 0 \\
\end{equation}
\noindent where $\int_{E}$ denotes integration over the subdomain $E \in \Omega_0$, and $\int_{e}$ denotes integration over the subdomain's boundary $e \in \Gamma_t$, which belongs to the traction boundary.

\subsection{Derivation of FPM trial and test functions} \label{ss:fpm_functions}
In each subdomain, trial and test functions are constructed by defining local discontinuous polynomial displacement vectors $\bm{u}^h$ and $\bm{v}$. Following the formulation in \cite{yang2019elementarily} and assuming that $\Omega_0$ is defined in $\mathbb{R}^{2}$, the trial function $\bm{u}^h \equiv \bm{u}^h(x,y)$ on the subdomain $E_0$, which encloses the point $P_0$, can be obtained by:
\begin{equation} \label{eq:fpm_trialfunc}
    \bm{u}^h(x,y) =
        \begin{bmatrix}
            u^h_x \\
            u^h_y
        \end{bmatrix} = 
            \begin{bmatrix}
                u_x^0 + \frac{\partial u_x}{\partial x} \bigg\rvert_{P_0} (x-x_0) + \frac{\partial u_x}{\partial y} \bigg\rvert_{P_0} (y-y_0)\\
                u_y^0 + \frac{\partial u_y}{\partial x} \bigg\rvert_{P_0} (x-x_0) + \frac{\partial u_y}{\partial y} \bigg\rvert_{P_0} (y-y_0)
            \end{bmatrix} (x, y) \in E_0
\end{equation}
\noindent where $(x_0, y_0)$ denote the coordinates of $P_0$,  $[u^0_x \quad u^0_y]^T$, the value of $\bm{u}^h$ at $P_0$, and $\left [\frac{\partial u_x}{\partial x} \quad \frac{\partial u_x}{\partial y} \quad \frac{\partial u_y}{\partial x} \quad \frac{\partial u_y}{\partial y} \right]^T \bigg\rvert_{P_0}$, the unknown derivatives of which are computed using the generalized Finite Difference Method \cite{liszka1980finite}. The support domain for $P_0$ is defined as the group of points $P_1$, $P_2$, $\dots$, $P_m$ that are involved in the subdomains $E_1$, $E_2$, $\dots$, $E_m$, which share an interface with subdomain $E_0$ (see Figure \ref{fig:partition}). We define the weighted discrete $L^2$ norm ($\bm{J}$) given in matrix form by:
\begin{equation} \label{eq:weighted_norm}
    \bm{J} = (\bm{A} \bm{\alpha} + \bm{u}_0 - \bm{u}_m)^T \bm{W}(\bm{A} \bm{\alpha} + \bm{u}_0 - \bm{u}_m)
\end{equation}
\noindent where
\begin{equation*}
    \bm{A} = \begin{bmatrix}
                x_1-x_0 & y_1-y_0 & 0 & 0 \\
                0 & 0 & x_1-x_0 & y_1-y_0 \\
                x_2-x_0 & y_2-y_0 & 0 & 0 \\
                0 & 0 & x_2-x_0 & y_2-y_0 \\
                \vdots & \vdots & \vdots & \vdots \\
                x_m-x_0 & y_m-y_0 & 0 & 0 \\
                0 & 0 & x_m-x_0 & y_m-y_0 \\
             \end{bmatrix}
\end{equation*}
\begin{equation*}
   \bm{a} = \left [\frac{\partial u_x}{\partial x} \quad \frac{\partial u_x}{\partial y} \quad \frac{\partial u_y}{\partial x} \quad \frac{\partial u_y}{\partial y} \right]^T \bigg\rvert_{P_0}
\end{equation*}
\begin{equation*}
    \bm{u}_0 =  \left [u_x^0 \quad u_y^0 \quad u_x^0 \quad u_y^0 \quad \dots \quad u_x^0 \quad u_y^0 \right]^T_{1\times 2m}
\end{equation*}
\begin{equation*}
    \bm{u}_m =  \left [u_x^1 \quad u_y^1 \quad u_x^2 \quad u_y^2 \quad \dots \quad u_x^m \quad u_y^m \right]^T
\end{equation*}
\begin{equation*}
    \bm{W} =  \begin{bmatrix}
            w_x^1 & 0 & 0 & \dots & 0 \\
            0 & w_y^1 & \ddots & \ddots & \vdots \\
            0 & \ddots & \ddots & \ddots & 0 \\
            \vdots & \ddots & \ddots & w_x^m & 0 \\
            0 & \dots & 0 & 0 & w_y^m
            \end{bmatrix}   
\end{equation*}
\noindent In the above matrices, $(x_i, y_i)$ denote the coordinates of point $P_i$, $\left[ u_x^i \quad u_y^i \right]^T$ and the value of $\bm{u}^h$ at $P_i$, $\left[ w_x^i \quad w_y^i \right]^T$ is the value of the weight function at $P_i$, $i = 1, 2, \dots, m$. The derivatives vector $\bm{a}$ is derived by solving the stationarity of $\bm{J}$ (Equation \ref{eq:weighted_norm}) to obtain:
\begin{equation} \label{eq:derivative_vector}
    \bm{a} = (\bm{A}^T \bm{W} \bm{A})^{-1} \bm{A}^T \bm{W} (\bm{u}^m - \bm{u}^0)
\end{equation}
\noindent By expressing $\bm{u}^m - \bm{u}^0$ as:
\begin{equation}
    \bm{u}^m - \bm{u}^0 = [\bm{I}_1 \quad \bm{I}_2] \bm{u}_E
\end{equation}
\noindent where 
\begin{equation}
\begin{split}
   \bm{u}_E = \left[ u_x^0 \quad u_y^0 \quad u_x^1 \quad u_y^1 \quad \dots \quad u_x^m \quad u_y^m \right]^T\\ \\
   \bm{I}_1 =
        \begin{bmatrix}
          -1 & 0\\
          0 & -1\\
          -1 & 0\\
          0 & -1\\
          \vdots & \vdots\\
          -1 & 0\\
          0 & -1
        \end{bmatrix}_{2m \times 2} \quad \bm{I}_2 =
        \begin{bmatrix}
         1 & 0 & \ldots & 0 \\
          0 & 1 & \ddots & \vdots \\
         \vdots & \vdots & \ddots & 0 \\
          0 & \ldots & 0 & 1
        \end{bmatrix}_{2m \times 2m}
\end{split}
\end{equation}
\noindent we can rewrite vector $\bm{a}$ as:
\begin{equation}
    \bm{a} = \bm{C} \bm{u}_E
\end{equation}
\noindent where 
\begin{equation} \label{eq:Cmat}
    \bm{C} = (\bm{A}^T \bm{W} \bm{A})^{-1} \bm{A}^T \bm{W} [\bm{I}_1 \quad \bm{I}_2]
\end{equation}
\noindent Finally, by substituting Equation (\ref{eq:Cmat}) in Equation (\ref{eq:fpm_trialfunc}), we obtain:
\begin{equation}
    \bm{u}^h = \bm{N} \bm{u}_E
\end{equation}
\noindent where the matrix $\bm{N}$ is the shape function of $\bm{u}^h$ in $E_0$ and is given by:
\begin{equation}
\begin{split}
    \bm{N} &= 
        \begin{bmatrix}
            x-x_0 & y-y_0 & 0 & 0 \\
            0 & 0 & x-x_0 & y-y_0 \\
        \end{bmatrix} \bm{C} + \bm{I}_3 \\
        \bm{I}_3 &= \begin{bmatrix}
            1 & 0 & 0 & \dots & 0 \\
            0 & 1 & 0 & \dots & 0
        \end{bmatrix}_{2\times(2m+2)}
\end{split}
\end{equation}

\noindent Since the presented FPM is based on the Galerkin weak form, the test function $v$ is derived in a similar manner to the trial function $u_h$. As there are no continuity requirements at the interfaces of the subdomains, the trial and test functions are discontinuous point-based simple and local polynomials. Due to the discontinuity, the Galerkin weak form leads to inconsistent and inaccurate solutions if continuity is not restored. It should be noted that instead of deriving the unknown derivatives by using the generalized Finite Difference Method, one may also employ the Differential Quadrature Method as was used in \cite{guan2021new}.

\subsection{Interior penalty numerical flux correction}
To solve the inconsistency problem, an interior penalty numerical flux correction is introduced, which is similar to discontinuous Galerkin FEM. The interior penalty numerical flux correction is commonly used in discontinuous Galerkin FEM to retrieve the continuity at the element interfaces \cite{arnold2002unified}. Similarly, consistency and accuracy in FPM are ensured by applying the interior penalty numerical flux correction to modify the Galerkin weak form in Equation (\ref{eq:galerkin_form}). This treatment can be viewed as applying a corrective internal force which acts on the interior interface $e \in \Gamma_h$ of subdomain $E$ and its neighbour, where $\Gamma_h$ denotes the interior part of $\Gamma$. Therefore, Equation (\ref{eq:galerkin_form}) is modified as:
\begin{equation} \label{eq:consistent_galerkin}
    \sum_{E \in \Omega_0} \int_{E} \bm{S} : \delta \bm{E} d\Omega - \sum_{E \in \Omega_0} \int_{E} \rho_0 (\bm{b} - \bm{\ddot u}^h) \cdot \bm{v} d\Omega - \sum_{e \in \Gamma_t} \int_{e} \bm{\Bar t} \cdot \bm{v} d\Gamma - \sum_{e \in \Gamma_h} \int_{e}{\bm{t}^{*} \cdot \llbracket \bm{v} \rrbracket d\Gamma} = 0 \\
\end{equation}
\noindent The numerical flux correction is defined using the incomplete interior penalty Galerkin (IIPG) method as in \cite{wang2022fragile}, which has an explicit physical meaning as the traction acting on the interior interface. The interior interface traction $\bm{t}^{*}$ is obtained by:
\begin{equation} \label{eq:correction_traction}
\begin{split}
    \bm{t}^{*} &= \{\ \bm{t}_h \} - \bm{\beta} \llbracket \bm{u} \rrbracket \\
    \bm{t}_h &= \frac{1}{2} \left( \bm{P}^{+} \bm{n}^{+}_{0} + \bm{P}^{-} \bm{n}^{-}_{0} \right) - \bm{\beta} \llbracket \bm{u} \rrbracket
\end{split}
\end{equation}
\noindent where $\llbracket \; \rrbracket$ and $\{\;\}$ are the jump and average operators, respectively. For any interior interface $e \in \partial E^+ \bigcap \partial E^-$ shared by the neighbouring subdomains $E^+$ and $E^-$, the operators act on an arbitrary quantity $w$ as:
\begin{equation} \label{eq:jump_avg_operators}
  \llbracket w \rrbracket = w |_e^{E^+} - w |_e^{E^-}, \quad
    \{w\} = \frac{1}{2} \left(w |_e^{E^+} + w |_e^{E^-}\right)
\end{equation}
\noindent In Equation (\ref{eq:correction_traction}), $\bm{\beta}$ is a second-order penalty tensor used to weakly enforce displacement continuity across the interior interfaces. It is defined as a diagonal matrix:
\begin{equation}
   \bm{\beta} = \frac{pE}{h_s} \bm{I} 
\end{equation}
\noindent where $p$ denotes a penalty coefficient, $E$ the Young’s modulus, and $h_s$ the subodmain characteristic length. $\bm{P}^{+}$, $\bm{P}^{-}$ denote the first Piola-Kirchhoff stress tensors for the two neighboring subdomains, while $\bm{n}_{0}^{+}$, $\bm{n}_{0}^{-}$ are the outward normal vectors to the interface of the neighboring subdomains at the reference configuration. The penalty tensor $\bm{\beta}$ is viewed as the interfacial stiffness. Therefore, the term $\bm{\beta} \llbracket \bm{u} \rrbracket$ in Equation (\ref{eq:correction_traction}) controls the contribution of separation between the neighboring subdomains, and $\bm{t}^{*}$ can be viewed as the exact traction acting on each interior interface, considering both the stress in the subdomains and the separation across the internal interface.

\subsection{Explicit integration in time}
The consistent Galerkin weak form in Equation (\ref{eq:consistent_galerkin}) can be written in matrix form as: 
\begin{equation} \label{eq:total_lagrangian}
    \bm{M} \ddot{\bm{u}} + \bm{F^{int}} = \bm{F^{ext}}
\end{equation}
\noindent where $\bm{u}$ denotes the displacements vector, $\bm{M}$ the mass matrix, $\bm{F^{int}}$ the internal forces, and $\bm{F^{ext}}$ the external forces (including body and traction forces). Applying mass lumping to diagonalize $\bm{M}$ allows solving Equation (\ref{eq:total_lagrangian}) explicitly using the central difference time integration scheme. Therefore, at each integration time step $k$, we obtain: 
\begin{equation} \label{eq:explicit_disp}
    \bm{u}_{k+1} = \bm{u}_{k} + dt_{k+1}\dot{\bm{u}}_k + \frac{1}{2dt^2_{k+1}} \ddot{\bm{u}}_k,
\end{equation}
\begin{equation}
    \dot{\bm{u}}_{k+1} = \dot{\bm{u}}_{k} + \frac{1}{2dt^2_{k+1}} (\ddot{\bm{u}}_{k+1} + \ddot{\bm{u}}_k),
\end{equation}
\begin{equation} \label{eq:explicit_stiffness_mod}
    \bm{F^{int}} = \sum_{E \in \Omega} \bm{F}_k^E - \sum_{e \in \Gamma_h} \Tilde{\bm{F}}_k^e = \sum_{E \in \Omega} \int_{E} \bm{B}^T \bm{S}_k d\Omega - \sum_{e \in \Gamma_h} \int_{e}{\llbracket \bm{N} \rrbracket ^{T} \bm{t}_{k}^{*} d\Gamma}, 
\end{equation}
\begin{equation} \label{eq:explicit_motion}
    \frac{1}{dt^2} \bm{M} \bm{u}_{k+1} = \bm{R}_k - \sum_{E \in \Omega} \bm{F}_k^E - \frac{1}{dt^2} \bm{M} (\bm{u}_{k-1} - 2\bm{u}_k)
\end{equation}
\noindent where $\bm{B}_k$ is the strain–displacement matrix at step $k$. It should be noted that there is no need to assemble the stiffness matrix $\bm{K}$ as Equations (\ref{eq:explicit_disp})--(\ref{eq:explicit_motion}) are solved explicitly. Internal forces are computed at the level of the non-overlapping subdomains. However, since explicit methods are only conditionally stable, an adequately small integration time step $dt$ must be selected to ensure the stability of the solution. An estimation of the critical stable time step is derived as in \cite{joldes2012stable}:
\begin{equation} \label{eq:time_step}
    \bar{dt}_{crit} = min_j \left (\frac{2}{\sqrt{\lambda^{max}_j}} \right)
\end{equation}
\noindent where $\lambda^{max}_j$ is the maximum eigenvalue of the stiffness matrix. The upper bound at a node $j$ enclosed by subdomain $E_j$ is calculated by:
\begin{equation}
    \lambda^{max}_j \leq m_j c^2 \bm{C}\bm{C}^T 
\end{equation}
\noindent where $m_j$ is the number of neighbouring points in the support domain of point $j$ and $c^2= \frac{\lambda + 2 \mu}{\rho_0}$ is the dilatation wave speed squared with $\lambda$. In this instance, $\lambda$, $\mu$ are the Lam\'e parameters and $\bm{C}$ the matrix of the FPM derivatives as given in Equation (\ref{eq:Cmat}).
\noindent It should be noted that the numerical flux correction term affects the numerical high-frequency eigenmodes \cite{noels2008explicit}. Therefore, the critical time step stability criterion must be modified. As described in \cite{noels2008explicit}, reducing the critical stable time step by the square root of the penalty coefficient is required. Therefore, Equation (\ref{eq:time_step}) should be updated to:
\begin{equation} \label{eq:new_time_step}
    dt_{crit} = \frac{1}{\sqrt{p}} \bar{dt}_{crit}, \quad p>0
\end{equation}

\section{Validation case studies} \label{sec:examples}

In this section, the explicit total Lagrangian FPM was employed to solve several validation case studies encompassing both 2D and 3D domains. The proposed method was evaluated in simulations of finite deformation at different levels of extension and compression. In all cases, hyperelastic material behaviour was considered using the standard neo-Hookean constitutive model with near-incompressibility. The strain energy density of this model was given in terms of the reduced invariants by:
\begin{equation} \label{eq:neohookean}
    W(J_1,J_3) = \frac{\mu}{2}(J_1-3) + \frac{K}{2}(J_3-1)^2
\end{equation}
\noindent where $\mu$ is the material's shear modulus, $K$ denotes its bulk modulus, and $J_1$, $J_3$ denote the first and third reduced invariants. The solutions obtained from FPM were compared with solutions obtained from FEM using nodally-integrated simplicial elements to avoid volumetric locking \cite{gee2009uniform}. To ensure a one-to-one comparison, the nodes of the FEM meshes were used for the FPM discretization. The discretization subdomains were generated using a dual polyhedral mesh generation algorithm \cite{garimella2014polyhedral,kim2014efficient}. Comparison was performed measuring the normalized root mean square error ($NRMSE$) of the obtained displacements from both FPM and FEM, which is given by:
\begin{equation}
    \textrm{NRMSE} = \frac{\sqrt{\displaystyle \frac{\displaystyle \sum_{j=1}^n \left( u^{h}_{j} - u^{ref}_{j} \right)^2}{n}}}{\displaystyle \max_{j} u^{ref}_{j}  - \displaystyle \min_{j} u^{ref}_{j}}
\end{equation}
\noindent where $u^h$ denotes the numerical solution obtained either by FPM or FEM, $u^{ref}$ the reference solution, and $n$ the number of nodes of the discretization. All simulations were performed on a laptop with Intel\textsuperscript{\textregistered} Core\textsuperscript{\texttrademark} i9-12900H CPU and 32GB RAM.

\subsection{Penalty coefficient effect on interior penalty numerical flux correction} \label{ss:penalty_effect}
\begin{figure}[hbt]
\captionsetup[subfigure]{justification=centering}
\begin{subfigure}[t]{0.475\textwidth}
    \includegraphics[width=\textwidth]{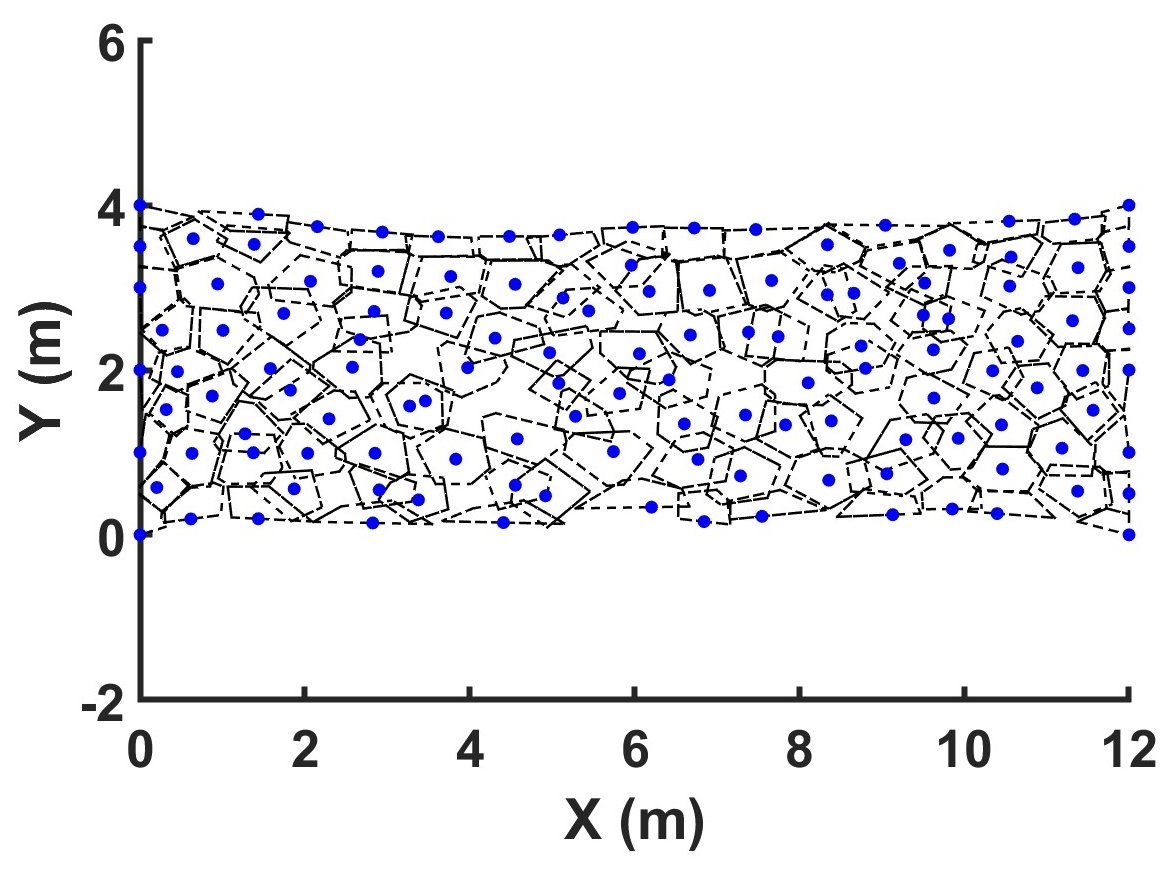}
    \caption{20$\%$ extension}
\end{subfigure}
\hfill
\begin{subfigure}[t]{0.475\textwidth}
    \includegraphics[width=\textwidth]{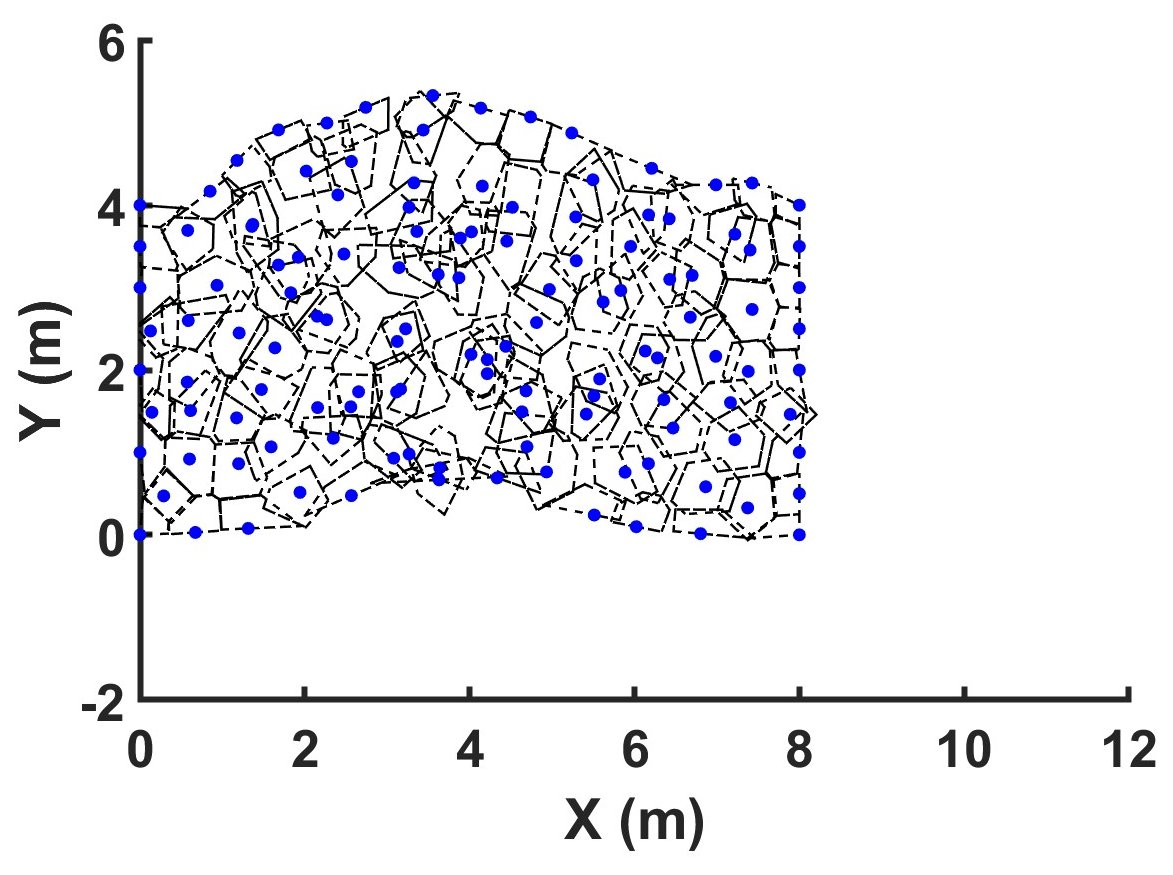}
    \caption{20$\%$ compression}
\end{subfigure}
\caption{Deformation of a 2D hyperelastic material using the explicit total Lagrangian Fragile Points Method (FPM) without interior penalty numerical flux correction ($p=0$). The inconsistency of the solution is evident due to the discontinuity of the FPM trial and test functions.}
\label{fig:inconsistent2D}
\end{figure}

In this case study, the effect of the penalty coefficient ($p$) on the interior penalty numerical flux correction was evaluated. Two simulation scenarios of a 2D neo-Hookean material with dimensions $10 \times 4$ m (nodes: 125) undergoing 20\% extension and 20\% compression were considered. Material density $(\rho) = 1000$ kg/m$^3$, Young's modulus $(\textrm{YM}) = 3$ kPa, and Poisson's ratio $(\nu) = 0.45$ were set. The material was constrained at $x = 0$ m applying zero displacement conditions $u_x = 0$ m, $u_y = 0$ m. At $x = 10$ m, the displacement condition $u_x = 2$ m was applied for the case of constrained extension, while $u_x = -2$ m was applied for constrained compression. Simulations were performed using $p = \{0, 10, 20, 50, 100\}$, where the obtained solution for $p=100$ was considered as the reference solution. When $p=0$, the solution was discontinuous at the subdomain interfaces, leading to the inconsistency problem observed in Figure \ref{fig:inconsistent2D}.

\begin{figure}[hbt]
\captionsetup[subfigure]{justification=centering}
\begin{subfigure}[t]{0.27\textwidth}
    \includegraphics[width=\textwidth]{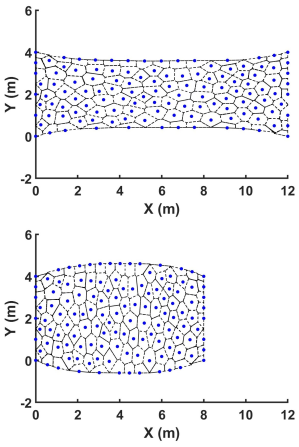}
    \caption{$p = 20$}
\end{subfigure}
\hfill
\begin{subfigure}[t]{0.27\textwidth}
    \includegraphics[width=\textwidth]{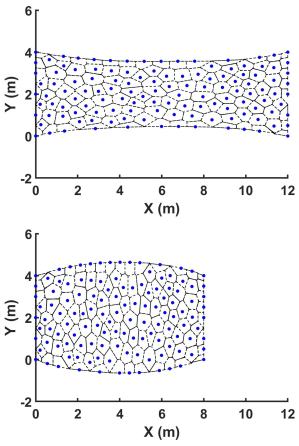}
    \caption{$p = 100$}
\end{subfigure}
\hfill
\begin{subfigure}[t]{0.45\textwidth}
    \includegraphics[width=\textwidth]{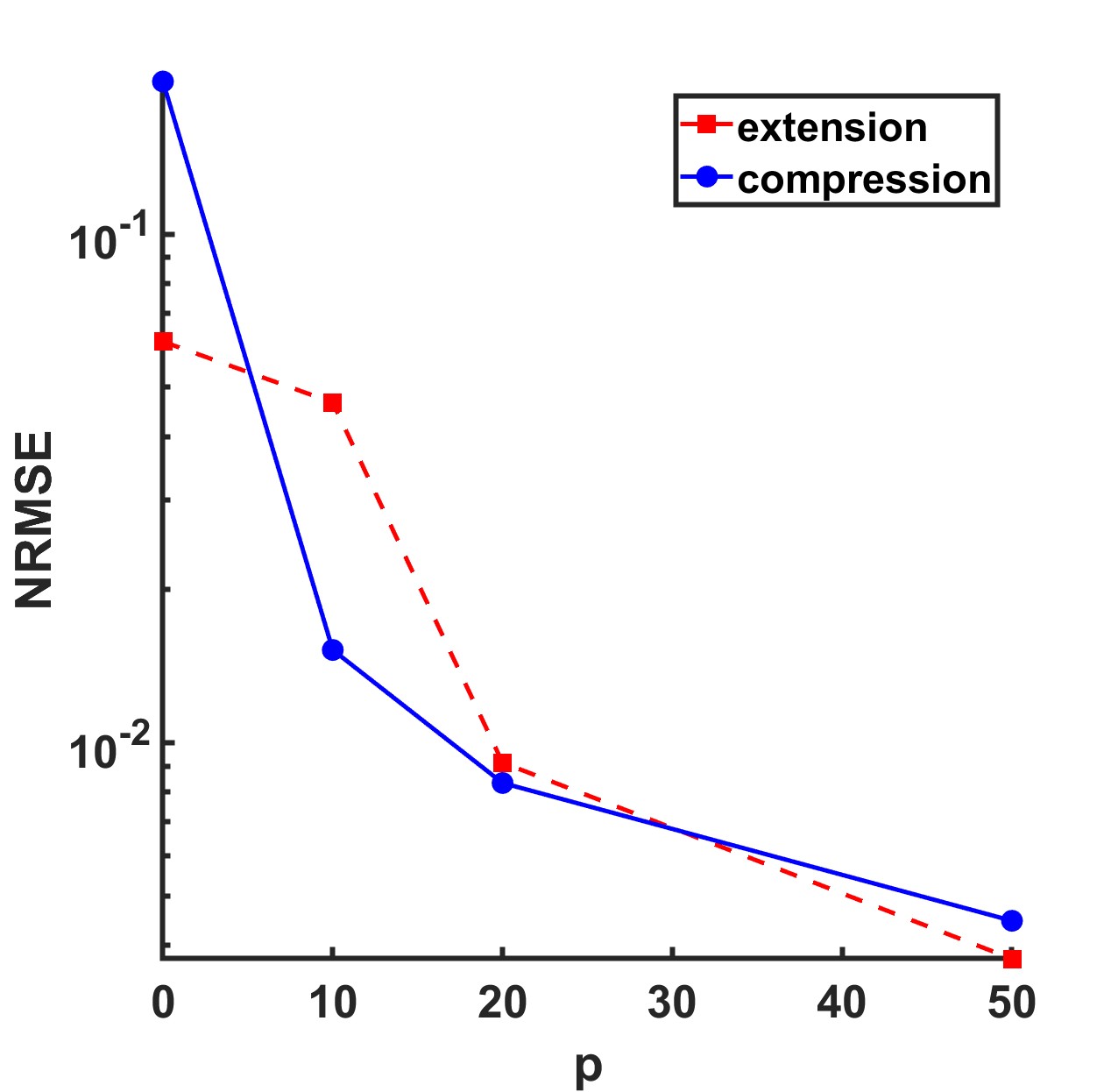}
    \caption{}
\end{subfigure}
\caption{Deformation of a 2D hyperelastic material at 20$\%$ extension (top) and 20$\%$ compression (bottom) using the explicit total Lagrangian FPM with interior penalty numerical flux correction with penalty coefficient (a) $p = 20$ and (b) $p = 100$. (c) $NRMSE$ for solutions with $p=0-50$ compared to a solution with $p=100$ for 20\% extension (red square) and 20\% compression (blue bullet).}
\label{fig:penalty2D}
\end{figure}

By increasing $p$, the discontinuity of the solution between the subdomains was reduced. For $p>0$, the consistency and accuracy were restored as shown in Figure \ref{fig:penalty2D}. From the NRMSE evaluation of the solutions with $p=0-50$ compared to the solution with $p=100$, accurate results were obtained for $p \geq 20$. For the case of constrained extension, NRMSE was found in the range $\left[6.48\mathrm{e}{-5}, 1.54\mathrm{e}{-2}\right]$, while for constrained compression, the NRMSE was found in the range $\left[5.32\mathrm{e}{-5}, 7.01\mathrm{e}{-2} \right]$. Recall from Equation (\ref{eq:new_time_step}) that as $p$ increases, the critical stable time step is reduced by $\sqrt{p}$. In order to ensure good accuracy and efficiency in the following examples, simulations were performed using $p=20$ as it was found to be a good trade-off value between accuracy and computational efficiency. 

\subsection{Unconstrained compression of a 3D hyperelastic block} \label{ss:unconstrained_comp3D}
\begin{figure}[hbt]
\captionsetup[subfigure]{justification=centering}
\begin{subfigure}[t]{0.55\textwidth}
    \includegraphics[width=\textwidth]{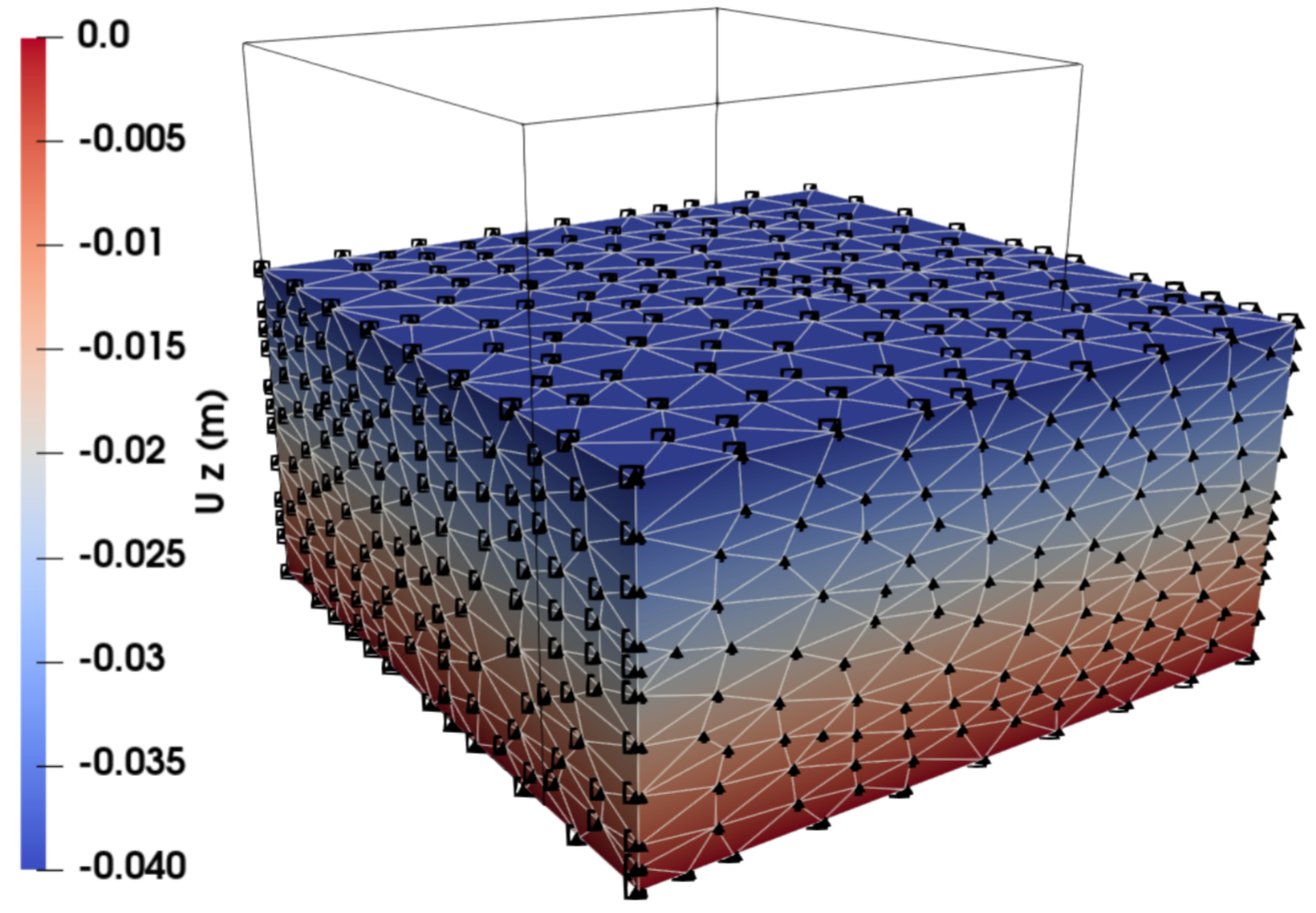}
    \caption{}
\end{subfigure}
\hfill
\begin{subfigure}[t]{0.4\textwidth}
    \includegraphics[width=\textwidth]{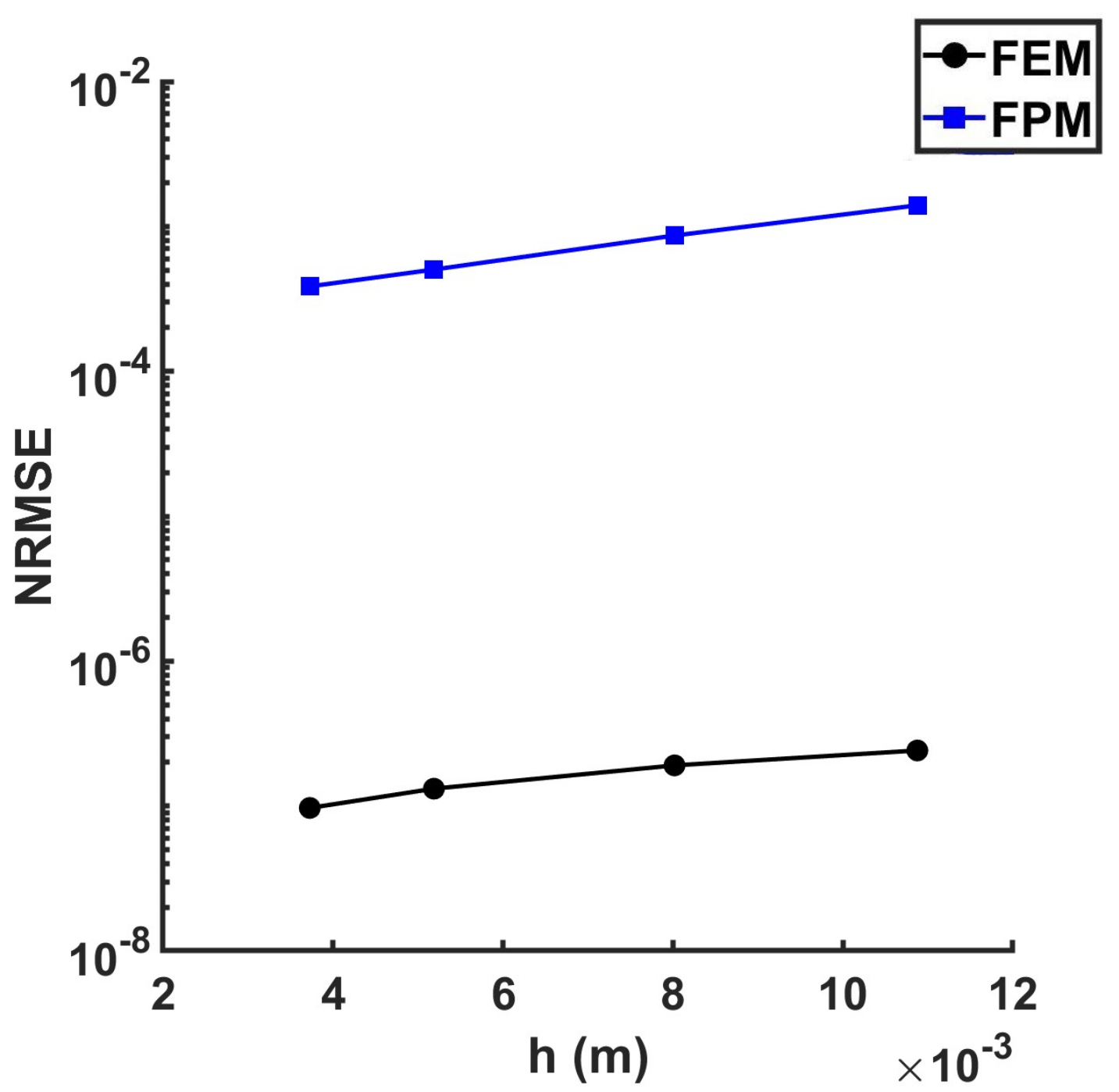}
    \caption{}
\end{subfigure}
\caption{(a) $40\%$ unconstrained compression of a 3D hyperelastic material using FPM (triangle) and FEM (square) for mesh with h = 1.09e-2 m. (b) $NRMSE$ of the displacement ($u_z$) obtained by FEM (black bullet) and FPM (blue square) for meshes with h = 1.09e-2 -- 3.73e-3 m compared to the analytical solution $u_z=-0.4z$, where $z$ denotes the z coordinate of the mesh nodes at the reference configuration.}
\label{fig:unconstrained_comp3D}
\end{figure}

In this case study, the unconstrained compression problem was solved for a 3D hyperelastic block with dimensions $0.1\times0.1\times0.1$ m$^3$. Unconstrained compression was simulated applying the following boundary conditions. $u_x = 0$ m at $x = 0$ m, $u_y = 0$ m at $y = 0$ m, $u_z = 0$ m at $z = 0$ m, and $u_z = 0.04$ m at $z = 0.1$ m. The block was modelled as a neo-Hookean material with $\rho = 1000$ kg/m$^3$, $\textrm{YM} = 3$ kPa, and $\nu = 0.45$. The problem was solved using FPM with $p=20$ and FEM for four different mesh resolutions with h = 1.09e-2 -- 3.73e-3, where h was the average nodal spacing (see Table \ref{tab:unconstrained_comp3D}). The NRMSE of the displacement was evaluated for both FPM and FEM solutions compared to the reference analytical solution $u_{ref} = -0.4z$, where $z$ denotes the Z-component of the nodes coordinates at the reference configuration. A qualitative comparison as well as the NRMSE convergence plots for FPM and FEM solutions are provided in Figure \ref{fig:unconstrained_comp3D}. The NRMSE for the different mesh resolution levels and the execution time for the computation of the FPM (t$_{FPM}$) and FEM (t$_{FEM}$) solutions are reported in Table \ref{tab:unconstrained_comp3D}. The execution time was normalized with respect to the execution time of the FEM solution for mesh resolution with $h=1.09e-2$ m (t$_{FEM}=2.84$ s).

\begin{center}
\begin{table*}[hbt]
\caption{$NRMSE$ and execution time report for Fragile Points Method (FPM) and Finite Element Method (FEM) solutions for the unconstrained compression of a 3D hyperelastic block.\label{tab:unconstrained_comp3D}}
\centering
\begin{tabular*}{500pt}{@{\extracolsep\fill}lccccc@{\extracolsep\fill}}
\toprule
\textbf{Nodes} & \textbf{h (m)} & \textbf{NRMSE$_{FPM}$} & \textbf{t$_{FPM}$}\tnote{$\dagger$} & \textbf{NRMSE$_{FEM}$} & \textbf{t$_{FEM}$}\tnote{$\dagger$}\\
\midrule
1639 & 1.09e-2 & 13.979e-4 & 1.39 & 2.416e-7 & 1\\
3843 & 8.02e-3 & 8.650e-4 & 3.81 & 1.907e-7 & 1.92\\
13144 & 5.19e-3 & 5.036e-4 & 19.93 & 1.317e-7 & 10.17\\
34064 & 3.73e-3 & 3.855e-4 & 83.2 & 0.969e-7 & 37.62\\
\bottomrule
\end{tabular*}
\begin{tablenotes}
\item[$\dagger$] FPM execution time (t$_{FPM}$) and FEM execution time (t$_{FEM}$) were normalized with respect to t$_{FEM}=2.84$ s, obtained from the FEM solution of the mesh resolution with h = 1.09e-2 m
\end{tablenotes}
\end{table*}
\end{center}

\subsection{Constrained extension of a 3D hyperelastic block} \label{ss:constrained_ext3D}
\begin{figure}[hbt]
\captionsetup[subfigure]{justification=centering}
\begin{subfigure}[t]{0.425\textwidth}
    \includegraphics[width=\textwidth]{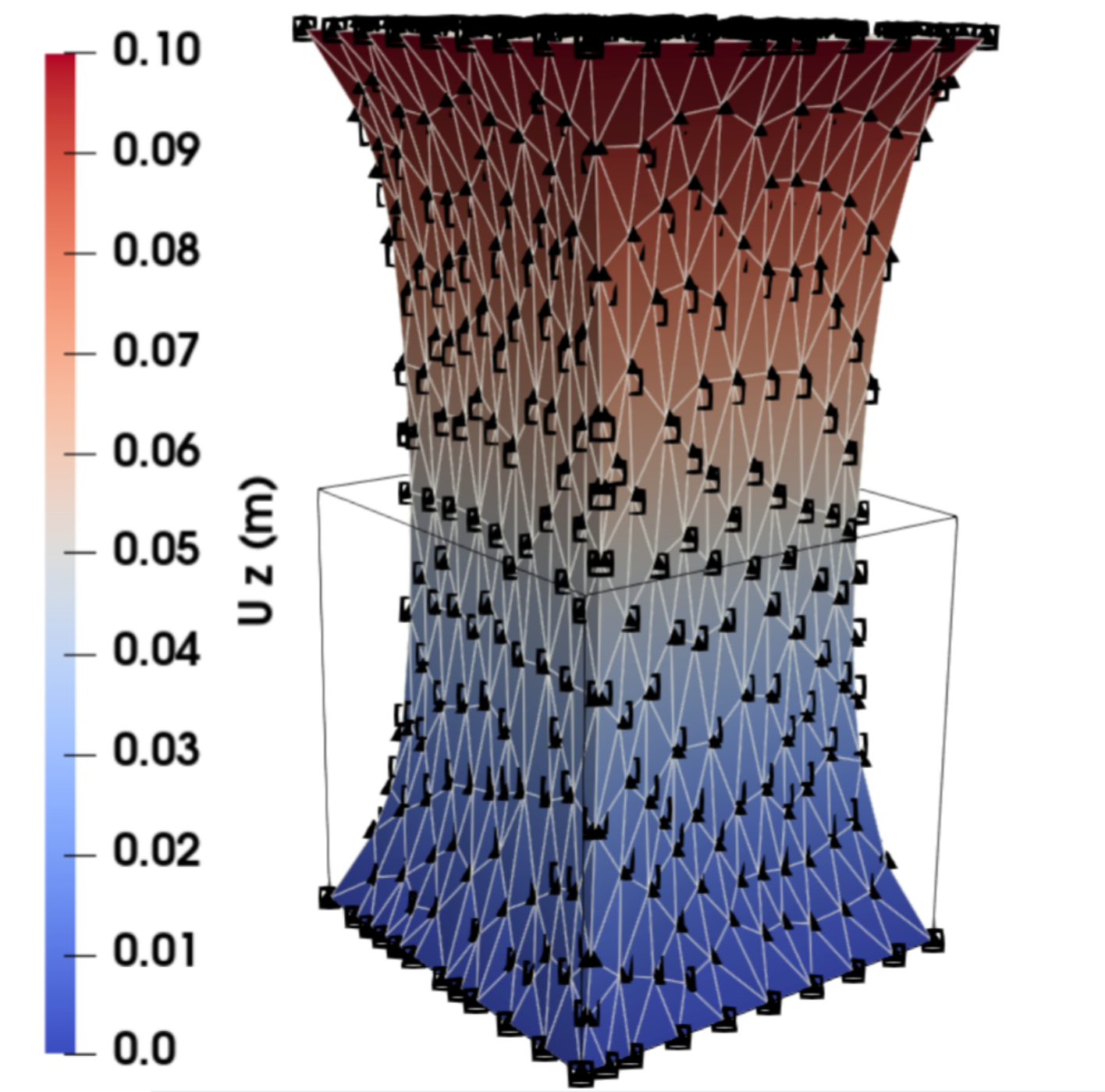}
    \caption{}
\end{subfigure}
\hfill
\begin{subfigure}[t]{0.475\textwidth}
    \includegraphics[width=\textwidth]{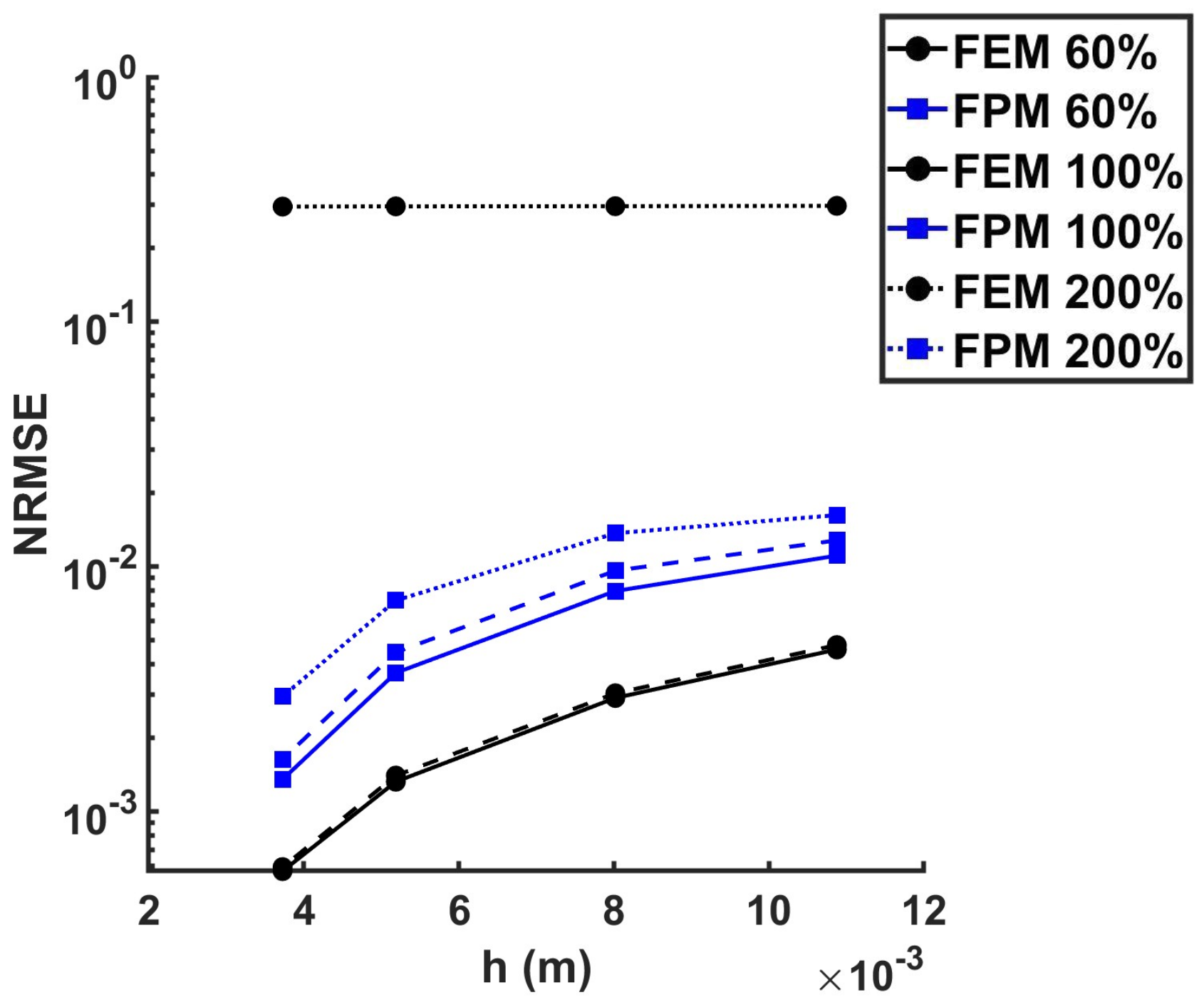}
    \caption{}
\end{subfigure}
\caption{(a) $100\%$ constrained extension of a 3D hyperelastic material using FPM (triangle) and FEM (square) for mesh with h = 1.09e-2 m. (b) $NRMSE$ of the displacement field ($\bm{u}$) obtained by FEM (black bullet) and FPM (blue square) for meshes with h = 1.09e-2 -- 3.73e-3 m compared to the solution for the reference mesh with h = 2.79e-3 m at $60\%$, $100\%$, and $200\%$ extension of the block's initial height.}
\label{fig:constrained_ext3D}
\end{figure}

In this case study, the constrained extension of the 3D neo-Hookean block in Example \ref{ss:unconstrained_comp3D} was considered. Zero displacement boundary conditions, $u_x = 0$ m, $u_y = 0$ m, $u_z = 0$ m, were applied at $z = 0$ m to fully constrain the bottom surface of the block. The top surface was constrained at $X-$ and $Y-$ directions applying $u_x = 0$ m and $u_y = 0$ m at $z = 0.1$ m. A fixed displacement $u_z = g$, where $g = {0.06, 0.1, 0.2}$ m was applied at $z=0.1$ m to simulate the extension of the block by $60\%$, $100\%$, and $200\%$ of its initial height, respectively. FPM solutions were obtained using $p=20$ and compared to solutions obtained from FEM for all the mesh refinement levels. Since an analytical solution was not available for this problem, the NRMSE for both FPM and FEM solutions was computed compared to a reference solution for a dense mesh with 80220 nodes and h = 2.79e-3 m. Prior to the NRMSE calculation, the FPM and FEM solutions for all the mesh refinement levels were mapped to the dense reference mesh using the Moving Least Squares approximation \cite{lancaster1981surfaces}. The NRMSE convergence plots are provided in Figure \ref{fig:constrained_ext3D}.

From the obtained NRMSE values for the different mesh refinement levels (see Table \ref{tab:constrained_ext3D}), the NRMSE convergence of FPM solution was in closer agreement with FEM compared to Example \ref{ss:unconstrained_comp3D}. This was mainly due to the increased mesh distortion in this example leading to numerical accuracy deterioration in the FEM solution. The mesh distortion effect in FEM accuracy was more evident in the case of the 200\% extension where the error was two orders of magnitude higher than the error in the 60\% and 100\% extensions. It should be noted that this accuracy deterioration was not observed in FPM where the error was maintained at the same order for all the extension cases.

\begin{center}
\begin{table*}[hbt]
\caption{$NRMSE$ report for FPM and FEM solutions for the constrained extension of a 3D hyperelastic block at 60\%, 100\%, and 200\% of its initial height.\label{tab:constrained_ext3D}}
\centering
\begin{tabular*}{500pt}{@{\extracolsep\fill}lcccccc@{\extracolsep\fill}}
\toprule
\multicolumn{1}{c}{\textbf{h (m)}} & \multicolumn{3}{c}{\textbf{NRMSE$_{FPM}$}} & \multicolumn{3}{c}{\textbf{NRMSE$_{FEM}$}}\\\cmidrule{2-4}\cmidrule{5-7}
 & 60\% & 100\% & 200\% & 60\% & 100\% & 200\% \\
\midrule
1.09e-2 & 1.110e-2 & 1.282e-2 & 1.625e-2 & 4.592e-3 & 4.788e-3 & 2.978e-1 \\
8.02e-3 & 7.956e-3 & 9.633e-3 & 1.370e-2 & 2.912e-3 & 3.043e-3 & 2.967e-1 \\
5.19e-3 & 3.688e-3 & 4.469e-3 & 7.282e-3 & 1.325e-3 & 1.403e-3 & 2.963e-1 \\
3.73e-3 & 1.354e-3 & 1.636e-3 & 2.954e-3 & 5.716e-4 & 5.935e-4 & 2.957e-1 \\
\bottomrule
\end{tabular*}
\end{table*}
\end{center}

\subsection{Constrained compression of a 3D hyperelastic block} \label{ss:constrained_comp3D}

\begin{figure}[hbt]
\captionsetup[subfigure]{justification=centering}
\begin{subfigure}[t]{0.55\textwidth}
    \includegraphics[width=\textwidth]{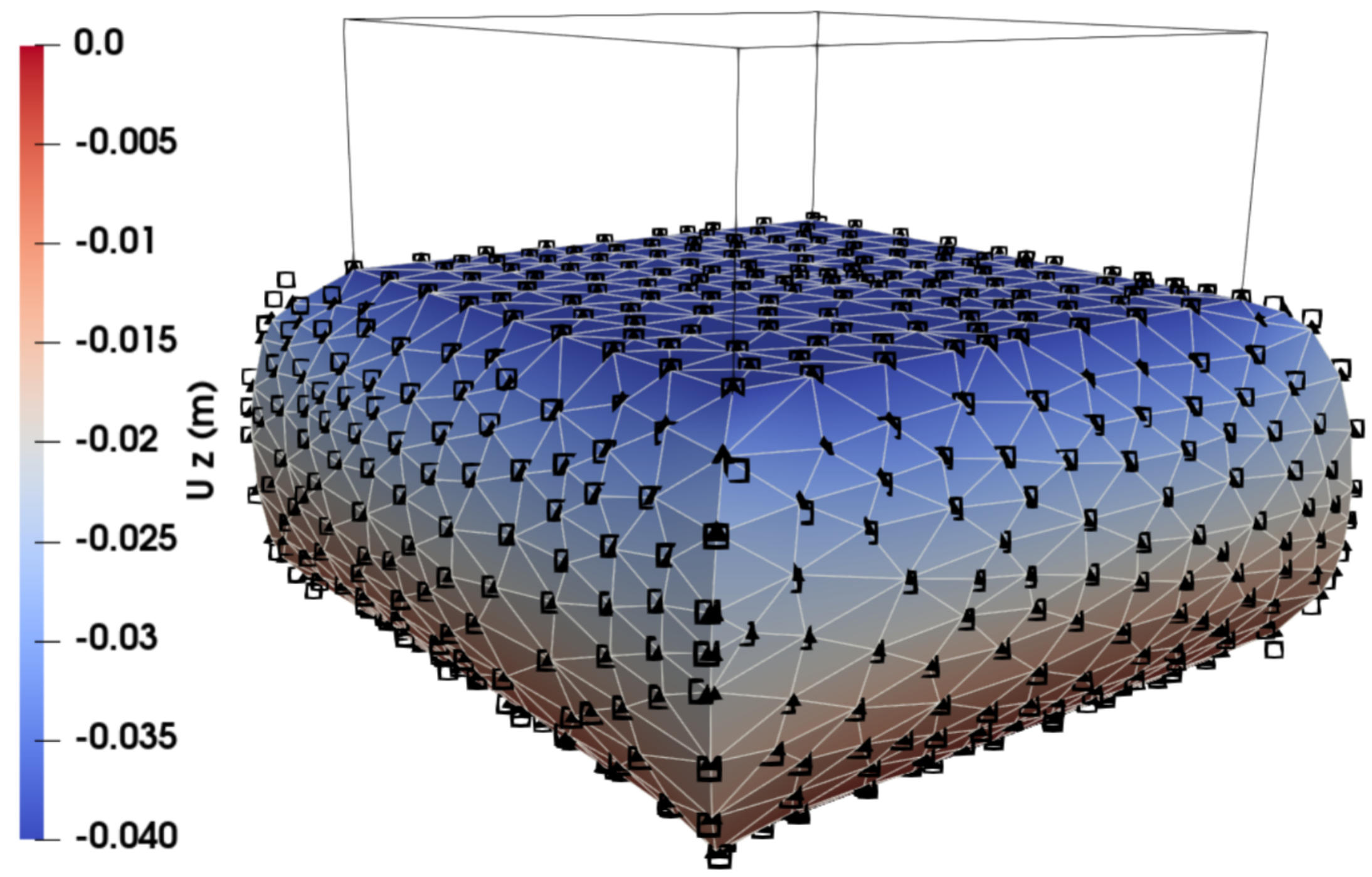}
    \caption{}
\end{subfigure}
\hfill
\begin{subfigure}[t]{0.42\textwidth}
    \includegraphics[width=\textwidth]{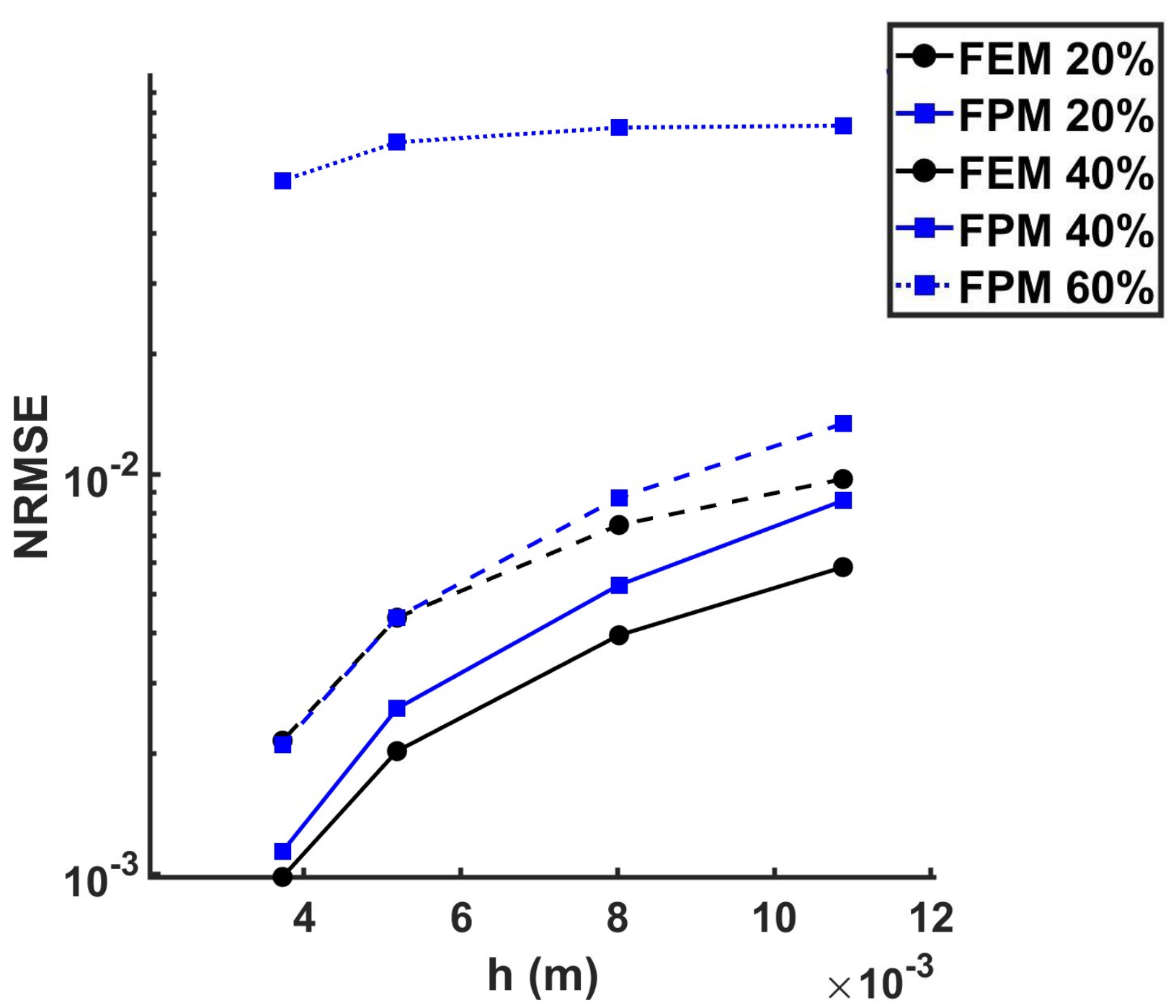}
    \caption{}
\end{subfigure}
\caption{(a) $40\%$ constrained compression of a 3D hyperelastic material using FPM (triangle) and FEM (square) for mesh with h = 1.09e-2 m. (b) $NRMSE$ of the displacement field ($\bm{u}$) obtained by FEM (black bullet) and FPM (blue square) for meshes with h = 1.09e-2 -- 3.73e-3 m compared to the solution for the reference mesh with h = 2.79e-3 m at $20\%$, $40\%$, and $60\%$ compression of the block's initial height. Note that FEM is unable to converge for compression larger than $40\%$.}
\label{fig:constrained_comp3D}
\end{figure}

In this case study, the setup of Example \ref{ss:constrained_ext3D} was used to simulate constrained compression of 20\%, 40\%, and 60\% of the initial height by applying the fixed displacement boundary condition $u_z = -l$, where $l = {0.02, 0.04, 0.06}$ m and $z = 0.1$ m at the top surface. The same material properties and zero displacement boundary conditions as in Example \ref{ss:constrained_ext3D} were applied in this scenario. Similarly, NRMSE was computed for FPM and FEM solutions; convergence plots are provided in Figure \ref{fig:constrained_comp3D}.

NRMSE values for FPM and FEM solutions of the constrained compression problem are given in Table \ref{tab:constrained_comp3D}. It is shown that the error has the same order of magnitude for FPM and FEM solutions for the cases of 20\% and 40\%. However, as can be seen in Figure \ref{fig:constrained_comp3D}a, the FEM solution at 40\% compression leads to highly distorted elements, especially at the corners of the geometry (square markers out of the deformed geometry). For 60\% compression, FPM converges in a smooth deformation state. However, NRMSE increases by one order of magnitude compared to the 20\% and 40\% compression cases. However, it should be noted that FEM solutions could not be obtained for 60\% since FEM could not converge for a compression rate higher than 40\%.

\begin{center}
\begin{table*}[hbt]
\caption{$NRMSE$ report for FPM and FEM solutions for the constrained compression of a 3D hyperelastic block at 20\%, 40\%, and 60\% of its initial height. Note that FEM failed to converge for 60\% constrained compression.\label{tab:constrained_comp3D}}
\centering
\begin{tabular*}{500pt}{@{\extracolsep\fill}lcccccc@{\extracolsep\fill}}
\toprule
\multicolumn{1}{c}{\textbf{h (m)}} & \multicolumn{3}{c}{\textbf{NRMSE$_{FPM}$}} & \multicolumn{3}{c}{\textbf{NRMSE$_{FEM}$}}\\\cmidrule{2-4}\cmidrule{5-7}
 & 20\% & 40\% & 60\% & 20\% & 40\% & 60\% \\
\midrule
1.09e-2 & 8.590e-3 & 1.338e-2 & 7.445e-2 & 5.856e-3 & 9.723e-3 & - \\
8.02e-3 & 5.281e-3 & 8.741e-3 & 7.363e-2 & 3.955e-3 & 7.464e-3 & - \\
5.19e-3 & 2.594e-3 & 4.385e-3 & 6.764e-2 & 2.029e-3 & 4.377e-3 & - \\
3.73e-3 & 1.139e-3 & 2.102e-3 & 5.423e-2 & 9.832e-4 & 2.155e-3 & - \\
\bottomrule
\end{tabular*}
\end{table*}
\end{center}

\section{Discussion} \label{sec:discuss}
Meshless methods have demonstrated their suitability compared to FEM for the large deformation of hyperelastic materials since the latter suffers from accuracy deterioration due to mesh distortion \cite{hu2011meshless,malkus1978mixed}. However, common meshless methods are usually more computationally expensive than FEM and require additional treatment for the exact imposition of boundary conditions. In this work, the Fragile Points Method (FPM) was considered for large deformation simulations as it delivers the advantages of meshless methods while overcoming their limitations. FPM is a novel meshless method which employs simple polynomials, which are local and discontinuous, as trial and test functions. Compared to other meshless methods, it has the advantage of simple and exact imposition of boundary conditions as in FEM. Moreover, the approximation is performed on compact support domains with simple trial and test functions that allow for accurate results with low-order integration. Therefore, the computational overhead of other meshless methods is significantly reduced for FPM. However, due to the discontinuity of trial and test functions in FPM, it is necessary to apply an interior penalty numerical flux correction to obtain consistent and accurate results. Since explicit time integration was used for the formulation of the total Lagrangian algorithm, the critical stable time step had to be reduced by the square root of the penalty coefficient ($p$) to ensure stability. Therefore, the value of $p$ should be chosen carefully to obtain accurate results without significantly reducing efficiency. In Example \ref{ss:penalty_effect}, different values of $p=0-100$ were investigated and $p=20$ was found to be a good trade-off value between accuracy and efficiency. 

The computational time of FPM compared to FEM was evaluated in Example \ref{ss:unconstrained_comp3D}. It was found that $t_{FPM} = \gamma t_{FEM}$, with $\gamma \in [1.39,2.2]$ being the FPM computational overhead which was mainly due to the evaluation of the first Piola--Kirchhoff stress tensor twice at each interface ($\bm{P}^+$, $\bm{P}^-$) during the computation of the interior penalty numerical flux correction at each integration time step. It should be noted that FPM execution time was compared with FEM execution time for FEM simulations employing simplicial elements. It is expected that if higher order elements were to be used, the FPM computational overhead should be lower compared to FEM. Therefore, it can be claimed that FPM is an efficient meshless method for large deformation simulation of hyperelastic materials.

In terms of accuracy, the performed experiments verified our initial hypothesis that FPM is a promising alternative to FEM when large mesh distortion occurs during deformation. Computing the normalized root mean square error, it was observed in case study \ref{ss:unconstrained_comp3D} that FEM had superior convergence compared to FPM for the simple case of the unconstrained compression of a hyperelastic block. This was expected as the mesh distortion was not severe in this case. However, when cases with higher mesh distortion were considered, e.g. case studies \ref{ss:constrained_ext3D} and \ref{ss:constrained_comp3D}, the NRMSE convergence of FEM was significantly reduced by four orders of magnitude. In contrast, the same order of magnitude was maintained for FPM convergence in all case studies. Most importantly, the FPM convergence was maintained even for the extreme cases of 200\% extension and 60\% compression of the initial height of the hyperelastic block. However, FEM convergence was not deteriorated for the case of 200\% extension. Furthermore, at 60\% compression a solution by FEM was not available since FEM failed to converge for compression rates larger than 40\%.

The simple neo-Hookean constitutive model was used in each case study despite it not being suitable for very large deformations. Nevertheless, FPM was able to produce smooth solutions even for the extreme deformation cases at 200\% extension and 60\% compression. The obtained results from this study demonstrated the capacity of the FPM implementation of the explicit total Lagrangian algorithm to simulate large deformations even for cases where large mesh distortion was involved and FEM failed to provide an accurate solution.

\section{Conclusion} \label{sec:conclude}
In this work, the Fragile Points Method (FPM) was employed to derive the explicit total Lagrangian algorithm and simulate the deformation of hyperelastic materials undergoing large deformation. Validation case studies were performed to evaluate the method against the standard Finite Elements Method (FEM). The results revealed that when mesh distortion is involved, FEM accuracy is deteriorated as expected, but FPM retains its accuracy. Moreover, FPM has minimal computational overhead and it leads to smooth solutions even for extreme deformation scenarios where FEM fails to converge. Therefore, we can conclude that FPM is a suitable meshless alternative to FEM for finite deformation of hyperelastic materials, especially for the simulation of severe deformation.

In future work, the proposed algorithm will be used to simulate soft tissue deformation during cardiac surgery procedures (e.g., catheter ablation) and cardiac electromechanics. It is also expected the proposed algorithm will be evaluated with more sophisticated constitutive materials such as the Holzapfel-Gasser-Ogden \cite{holzapfel2000new} and Guccione \cite{guccione1991passive} constitutive models. In the context of soft tissue simulation where time restrictions of the corresponding clinical application may apply, the work presented in \cite{meister2020deep} is of great interest. The use of neural networks as function approximators was proposed to accelerate the explicit total Lagrangian algorithm by allowing the use of a time step up to 20 times larger than the critical stable time step. It is expected that by using the explicit total Lagrangian FPM for training data generation, the neural network could generate accurate results for a wide range of deformations.

\subsection*{Acknowledgements}
This work has received funding from the European Union’s Horizon 2020 research and innovation programme under the Marie Skłodowska-Curie grant agreement PhyNeTouch No 101024463.

\subsection*{Conflict of interest}
The authors declare no potential conflict of interests.

\bibliographystyle{unsrt}  
\bibliography{main}

\end{document}